\documentclass[11pt]{article}
\usepackage{hyperref}
\hypersetup{hypertex=red,
	colorlinks=true,
	linkcolor=red,
	anchorcolor=red,
	citecolor=red}

\usepackage{psfrag}
\usepackage{epstopdf}
\usepackage{epsfig}
\usepackage{verbatim}
\usepackage{fancyhdr}
\usepackage{subfigure}
\usepackage{indentfirst}
\usepackage{booktabs}
\usepackage{color}
\usepackage{mathrsfs}
\usepackage{graphics}
\usepackage{graphicx}
\usepackage{listings}
\usepackage{xcolor}
\usepackage{float}
\usepackage{multirow}

\usepackage{amsmath}
\usepackage{amssymb}
\usepackage{amscd}
\usepackage{amsthm}

\usepackage[all]{xy}

\usepackage{algorithm}  
\usepackage{algorithmic}  
\usepackage{enumerate}

\usepackage{xypic}  

\usepackage{DG_macro}
\allowdisplaybreaks

\begin{document}

	\title{\textbf{A sparse approximate inverse 
		for triangular matrices based on Jacobi iteration}}
%	
%	\author{ Zhongjie Lu 
%		\thanks{School of Mathematical Sciences, 
%			University of Science and Technology of China, Hefei, Anhui 230026, China.
%			Email: zhjlu@ustc.edu.cn}}

    \author{ Zhongjie Lu 
    		\thanks{Dipartimento di Matematica, Universit\`{a} di Pavia, via Ferrata 5, 27100, Pavia, Italy.
    			Current address: School of Mathematical Sciences, 
    			University of Science and Technology of China, Hefei, Anhui, 230026, China. 
    			Email: zhjlu@ustc.edu.cn}}

	\date	

	\maketitle

\begin{abstract}
In this paper, we propose a simple sparse approximate inverse 
for triangular matrices (SAIT).
Using the Jacobi iteration method,
we obtain an expression of the exact inverse of triangular matrix,
which is a finite series.
The SAIT is constructed based on this series.
We apply the SAIT matrices to iterative methods with ILU preconditioners.
The two triangular solvers in the ILU preconditioning procedure
are replaced by two matrix-vector multiplications,
which can be fine-grained parallelized.
We test this method by solving some linear systems and eigenvalue problems
with preconditioned iterative methods.
\end{abstract}


\section{Introduction}

The incomplete LU (ILU) factorization is a type of 
general-purpose preconditioning techniques
for sparse linear systems.
There are two main problems in its parallelization.
The first is the parallel generation of ILU factors.
Many theories and techniques have been used 
to improve its parallelization.
Here are some references on this problem
\cite{abdelfattah2016linear, anzt2018parilut,chow2015fine,
	henon2006parallel, hysom2001scalable,monga2001generalized}.
%
%
%In our previous paper \cite{iterilu},
%we propose an iterative ILU (IterILU) factorization in matrix form. % Algorithm \ref{IterILU}.
%The factors can be generated by several matrix-matrix multiplications
%and some dropping rules,
%which are of fine-grained parallel operations.
%The IterILU has similar preconditioning effect 
%to the conventional ILU factorizations in solving many problems.
%
%
The second problem is solving triangular systems 
in the preconditioning procedure.
It is easy to solve exactly a triangular system 
using forward or backward substitution method.
However, this is a highly sequential process,
and it may be executed many times in solving a system.
This is an obstacle to exploit the performance of 
a parallel computing platform sufficiently.
This is the main problem we shall study in this paper.

The simplest parallelization strategy in solving triangular systems
is to compute the summations in the substitution methods in parallel:
\begin{equation}\nonumber %label{}
	\begin{split}
		x_i = \frac{1}{a_{ii}} \left( b_i - \sum_{k=1}^{i-1}a_{ik}x_k \right)
		\qq \text{or} \qq
		x_i = \frac{1}{a_{ii}} \left( b_i - \sum_{k=i+1}^{n}a_{ik}x_k \right).
	\end{split}
\end{equation}
This is a low-level concurrent method,
especially for sparse triangular matrices.
Its corresponding blocking method can improve the performance.
The level-scheduling method with reordering techniques
improve the parallelism further \cite{anderson1989solving}.
These methods are exact solvers for triangular systems.

Another type of strategy is to replace the exact solutions of 
triangular systems by approximate ones.
This ideal is based on the fact that
the incomplete LU factorization is incomplete.
This means that even the exact solutions of the triangular systems in ILU
are not 'exact' and these exact solutions could tolerate some errors.
The Jacobi iteration method is an easy way 
to obtain approximate solutions \cite{chow2018using, chow2015fine},
and it mainly involves matrix-vector multiplications and vector additions,
which are of fine-grained parallelism.
An alternative method is to use sparse approximate inverses (SAI),
which is based on the decay of inverses of sparse matrices 
\cite{demko1984decay,eijkhout1988decay, nabben1999decay}.
With SAIs,
the two triangular solvers in ILU preconditioning 
can be replaced by 
two matrix-vector multiplications.
There are several ways to compute the SAI of a matrix,
for example Frobenius norm minimization and incomplete inverse factorization
\cite{benzi2002preconditioning, benzi1996sparse,  benzi1999comparative, saad2003iterative}.
Also, there are some SAIs specially designed for triangular systems
\cite{anzt2018incomplete,janna2010block,van1999scalable}.
In \cite{ gustafsson1995completely,van1982vectorizable},
the truncated Neumann expansions play a similar role to the SAIs 
of ILU factors.
Some recent processes of this topic are introduced in
\cite{abdelfattah2016linear}.

When deriving the method we propose in this paper,
we compress the Jacobi iteration method to a series firstly.
As the iterative matrix in Jacobi iteration method for triangular systems 
is strictly triangular,
the series is finite.
We take the truncation of this series as the
approximate inverse of triangular matrix
and use some dropping strategies to make 
the inverse keep enough sparsity.
We rewrite the truncation of the series as
a recursive formula,
by which we can avoid computing and storing high-order matrix powers.
Then the SAIT algorithm is written in a quite brief way
and it mainly involves sparse matrix-matrix 
multiplications and some dropping rules.
The SAIT matrices can be applied to ILU preconditioning.
The two triangular solvers can be replaced by 
two matrix-vector multiplications with a pair of SAIT matrices,
which can be fine-grained parallelized.

This paper is organized as follows.
In section 2, we propose an exact inverse for triangular matrices
base on Jacobi iteration
and then construct the general formulation of the SAIT algorithm.
In section 3, we use two dropping strategies in the SAIT algorithm.
In section 4, we take the SAIT matrices as preconditioners 
in solving linear systems and eigenvalue problems with 
preconditioned iterative methods.
In section 5, there are some conclusions.

\section{The basic idea of the SAIT algorithm}
Let $ T\in \mbR^{n\times n} $ be a triangular matrix.
We solve the linear equation
\begin{equation}\nonumber %label{}
	\begin{split}
		Tx = b
	\end{split}
\end{equation}
by Jacobi iteration method
\begin{equation}\label{Jacobi_iter}
	\begin{split}
		x^k = D^{-1}(D-T)x^{k-1}+D^{-1}b.
	\end{split}
\end{equation}
Here, $ D $ is the diagonal matrix of $ T $.
%As $ T $ is triangular, 
%the iteration solution converges
%to the exact solution within at most $ n $ steps.
For the sake of simplicity,
we denote $ \tT = I - D^{-1}T $.
$ \tT $ is a strictly triangular matrix.
Then, \eqref{Jacobi_iter} becomes
\begin{equation}\nonumber %label{}
	\begin{split}
		x^k = \tT x^{k-1} + D^{-1}b.
	\end{split}
\end{equation}
Letting the initial datum $ x^0 = 0 $,
we list the solution of each iteration:
\begin{equation}\nonumber %label{}
	\begin{split}
		x^1 &= D^{-1}b\\
		x^2 &= \tT x^1 +D^{-1}b = (\tT +I)D^{-1}b\\
		x^3 &= \tT x^2 +D^{-1}b = (\tT^2 + \tT + I)D^{-1}b\\
		& \cdots\\
		x^k &= \tT x^{k-1} +D^{-1}b = \left( \sum_{i=0}^{k-1} \tT^i D^{-1}\right) b.\\
	\end{split}
\end{equation}
Denoting 
\begin{equation}\label{M_poly}
	\begin{split}
		M_k = \sum_{i=0}^{k-1} \tT^i D^{-1},
	\end{split}
\end{equation}
then we have 
\begin{equation}\nonumber %label{}
	\begin{split}
		x^k = M_k b.
	\end{split}
\end{equation}
As $ \tT $ is a strictly triangular matrix,
we have $ \tT^k = 0 $ when $ k\geq n $.
Then we know that
the exact inverse of $ T $ is a finite series of \eqref{M_poly},
i.e.
\begin{equation}\label{T_inverse}
	\begin{split}
		T^{-1} = M_n \equiv \sum_{i=0}^{n-1} \tT^i D^{-1}.
	\end{split}
\end{equation}

%In the ILU preconditioning,
%the triangular factors are not the exact factorization 
%of the original matrix.
%There are some errors contained in them.
%We denote the error of $ T $ by $ T_\epsilon $.
%Then the ideal solution should be 
%\begin{equation}\nonumber %label{}
%\begin{split}
%%
%x_{ideal} = (T+T_\epsilon)^{-1} b.
%%
%\end{split}
%\end{equation}
%Even though we solve the exact solution of the system $ Tx = b $ 
%($ x = T^{-1}b \equiv M_{n}b $),
%there is still a gap between $ x $ and $ x_{ideal} $.
%Based on this fact,
%if we impose an appropriate perturbation $ M_\epsilon $
%on the exact inverse of $ T $
%and get an approximate solution of $ x $ 
%by $ x_\epsilon = (T^{-1}+M_\epsilon)b \equiv (M_{n}+M_\epsilon)b $,
%the error magnitude of $ \nm{x_\epsilon - x_{ideal}} $
%can be similar to $ \nm{x - x_{ideal}} $.
%Then $ x $ can be replaced  by $ x_\epsilon $ 
%in the preconditioning procedure.
% 
%The perturbation $ M_\epsilon $ can be some dropping 
%and truncation rules.

In the following, 
we construct approximate inverses for the triangular matrix $ T $
based on its exact inverse \eqref{T_inverse}.
The approximation comes from two aspects, truncation and dropping.
Let we study a truncated term  $ M_m $ of \eqref{T_inverse},
where $ m<n-1 $ and
\begin{equation}\label{M_m}
	\begin{split}
		 M_m = D^{-1} +  \tT D^{-1} + \tT^2 D^{-1}+  \cdots + \tT^m D^{-1}.
	\end{split}
\end{equation}
In order to avoid computing and storing 
the high-order terms $ \tT^i $ $ (i=1,\cdots,m)$ in \eqref{M_m},
we rewrite $ M_k $ as an equivalent form:
\begin{equation}\label{Horner}
	\begin{split}
		M_0D &= I,\\
		M_1 D &= \tT +I,\\
		M_2 D &= \tT(\tT+I)+I,\\
		M_3 D &=  \tT(\tT(\tT+I)+I)+I,\\
		& \cdots\\
	\end{split}
\end{equation}
This is the Horner's method in calculating the values of polynomial functions.
The same method was also proposed by a Chinese mathematician 
\emph{Qin Jiushao}
%in his book \emph{Mathematical Treatise in Nine Sections} 
%(Chinese: $\ll$数书九章$\gg$)
in $13^{th}$ century Song dynasty \cite{Wu2000GrandSeries}.
Furthermore, we present \eqref{Horner} in a recursive formula:
\begin{equation}\label{Recursive}
	\begin{split}
		M_0D &= I,\\
		M_1 D &= \tT (M_0D) +I,\\
		M_2 D &= \tT (M_1D) +I,\\
		M_3 D &= \tT (M_2D) +I,\\
		& \cdots\\
		M_m D &= \tT (M_{m-1}D)+I.\\
	\end{split}
\end{equation}
To make the matrix more sparse, 
we perform dropping rules after 
each multiplication by $ \tT $ in \eqref{Recursive}.
%As $ M_k $ contains high-order term $ \tT^i $ $ (i=1,\cdots,m)$,
%too many nonzeros are contained in the finial result.
%Dropping some of them can make the storage at a reasonable level.
This strategy is summarized in Algorithm \ref{SAIT_basic}.
\begin{algorithm}  
	\caption{Sparse Approximate Inverse for Triangular matrices
		base on Jacobi iteration (SAIT)}  
	\begin{algorithmic}[1]  
		\STATE	Set the initial data $ M=I $,
		let $ D $ be the diagonal matrix of $ T $
		and let $ \tT = I-D^{-1}T $.
		\FOR { $ k = 1, 2, \cdots,m $}
		\STATE  $ M = \tT M + I $
		\STATE  execute dropping rules on  $ M $
		\ENDFOR
		\STATE $ M = MD^{-1}$  
	\end{algorithmic}  
	\label{SAIT_basic}
\end{algorithm}

\section{Dropping strategies}
We use two dropping strategies,
threshold-based and pattern-based, to concretize Algorithm \ref{SAIT_basic}.

\textbf{Threshold-based.}
In step 3 of Algorithm \ref{SAIT_basic},
we drop the entries whose magnitudes are small than 
$ \tau $ $ (\tau<1) $ in $ M $.
As the diagonal entries of $ M $ are all units,
the setting $ (\tau<1) $  can guarantee
the matrix $ M $ being always full-rank during the iterations.
This strategy is presented in Algorithm \ref{SAIT_Thr}.
We paste the Matlab code for this algorithm
in Appendix \ref{SAIT_Thr_code}.

\begin{algorithm}  
	\caption{Threshold-based Sparse Approximate Inverse for Triangular matrices
		SAIT\_Thr$ (\tau, m) $}  
	\begin{algorithmic}[1]  
		\STATE	Set the initial data $ M=I $,
		let $ D $ be the diagonal matrix of $ T $
		and let $ \tT = I-D^{-1}T $.		
		\FOR { $ k = 1, 2, \cdots,m $}
		\STATE  $ M = \tT M + I $
		\STATE  drop the entries in  $ M $ whose magnitudes are small than $ \tau $ $ (\tau<1) $
		\ENDFOR
		\STATE $ M = MD^{-1}$  
	\end{algorithmic}  
	\label{SAIT_Thr}
\end{algorithm}

\textbf{Pattern-based.}
We set a fixed sparsity pattern $ \vS $ in advance,
and drop the entries out of $ \vS $ in each iteration.
The patterns of the power functions of the original matrix
are often used in computing  sparse approximate inverses \cite{saad2003iterative}.
We study the sparsity pattern of $ T^p $ by
\begin{equation}\nonumber %label{}
	\begin{split}
		T^p = ((T-D) + D )^p = \sum_{i=0}^p C_p^i D^{p-i}(T-D)^i.
	\end{split}
\end{equation}
Since $ (T-D) $ shares the same pattern with $ \tT $,
compared with \eqref{M_poly},
we know that $ T^p $ has the same pattern with $ M_{p+1} $.
As $ M_p $ tends to the exact inverse with $ p \to n $,
this pattern can be regarded as a truncated pattern of the exact inverse.
When $ p=0 $, it is the Jacobi preconditioner,
while when $ p=1 $, it is the pattern of the original matrix $ T $.
The pattern of $ T^p $ can be obtained by Algorithm \ref{SAIT_basic} 
without dropping any entries in the first $ p $ iterations.
Then,  we keep the pattern of $ T^p $,
and in the following iterations,
drop the entries out of this pattern.
The number $ p $ should not be too large,
as the number of nonzeros in the pattern of $ M_p $ increases rapidly with $ p $ increasing.
This dropping strategy is summarized in Algorithm \ref{SAIT_Pat}.
\begin{algorithm}  
	\caption{Pattern-based Sparse Approximate Inverse for Triangular matrices
		SAIT\_Pat$ (p,m) $}  
	\begin{algorithmic}[1]  
		\STATE	Set the initial data $ M=I $,
		let $ D $ be the diagonal matrix of $ T $
		and let $ \tT = I-D^{-1}T $.		
		\FOR { $ k = 1, 2, \cdots,p $}
		\STATE  $ M = \tT M + I $
		\ENDFOR
		\STATE keep the sparsity pattern $ \vS $ of $ M $
		\FOR { $ k = 1, 2, \cdots,m $}
		\STATE  $ M = \tT M + I $
		\STATE  drop the entries in  $ M $ out of the pattern $ \vS $
		\ENDFOR
		\STATE $ M = MD^{-1}$  
	\end{algorithmic}  
	\label{SAIT_Pat}
\end{algorithm}

%\subsection{Hybrid strategy}
%\begin{algorithm}  
%	\caption{ Sparse Approximate Inverse for Triangular systems}  
%	\begin{algorithmic}[1]  
%		\STATE	Set the initial data $ M=I $,
%		 let $ D $ be the diagonal matrix of $ T $
%		 	and let $ \tT = I-D^{-1}T $.
%		\FOR { $ k = 1, 2, \cdots,m $}
%		\STATE  $ M = \tT M + I $
%		\ENDFOR
%		\FOR { $ k = 1, 2, \cdots,p $}
%		\STATE  $ M = \tT M + I $
%		\STATE  drop the entries in  $ M $ whose magnitudes are small than $ \tau $
%		\ENDFOR
%		
%		\STATE $ M = MD^{-1}$  
%	\end{algorithmic}  
%	\label{SAIT_hybrid}
%\end{algorithm} 

\subsection{Balance between the accuracy and iteration counts}

Let $ L $ and $ U $ be the two triangular factors 
of incomplete LU (ILU) factorization of a matrix $ A $,
and let $ M_L $ and $ M_U $
be the SAIT matrices for them respectively.
In ILU preconditioning,
the multiplication $ LU $ can be regarded as 
an approximation of $ A $, i.e. $ LU \approx A $.
If we take  $ M_L $ and $ M_U $ 
as the approximations of $ L^{-1} $ and $ U^{-1} $, respectively,
the multiplication $ M_U M_L $ can be take 
an approximate inverse of $ A $, i.e.
\begin{equation}\nonumber %label{}
	\begin{split}
		M_U M_L \approx U^{-1}L^{-1} \approx A^{-1}.
	\end{split}
\end{equation}
When the ILU is taken as a preconditioner in an iterative method,
we replace the two triangular equations $ U (L x) = y $
by two matrix-vector multiplications $ x = M_U(M_L y) $.

In the Algorithm \ref{SAIT_Thr} and \ref{SAIT_Pat},
if we use smaller threshold $ \tau $ and larger $ p $,
we can obtain more accurate approximate inverses,
since the SAIT matrices contains more nonzero entries.
However, higher accuracy of SAIT matrices do not mean
shorter runtime in solving a system with this preconditioner.

We list the following terms related to an iterative solver 
with the SAIT preconditioners:
\begin{itemize}
	\item  $ iter(M_L, M_U) $:
	the iteration count of an iterative solver 
	with $ M_L $  and $ M_U $ in preconditioning procedure.
	
	\item  $ comp(M_L, M_U) $: the runtime of computing the two matrix-vector multiplications $ M_U(M_L y) $.
	
	\item  $ othercomp $: the runtime of computing other terms in each iteration 
	except the preconditioning procedure.
\end{itemize}
Then, the runtime of a preconditioned iterative method 
can be estimated roughly by the following formula:
\begin{equation}\label{runtime}
	\begin{split}
		runtime = iter(M_L, M_U)\times\bigg(comp(M_L, M_U) + othercomp \bigg).
	\end{split}
\end{equation}
We define the ratios of a SAIT matrix $ M_T $ compared 
with its original matrix $ T $:
\begin{equation}\nonumber %label{}
	\begin{split}
		r  = \frac{nnz(M_T)}{nnz(T)},
	\end{split}
\end{equation}
where $ nnz(\cdot) $ denotes the number of nonzeros in a matrix.
This value is to illustrate the number of nonzeros in a SAIT matrix,
and it correlates to 
the runtime of the part $ comp(M_L, M_U) $.

When we use iterative methods to solving a linear system,
if the iterative method and the computer environment are fixed,
the term $ othercomp $ is usually fixed.
If we use more accurate SAIT matrices,
the iteration count $ iter(M_L, M_U) $ decreases.
However,
we can not expect that the SAIT preconditioners can reduce the iteration count
less than exact triangular solvers.
It is probable that after exceeding some point,
more nonzeros in $ M_L $ and $ M_U $
can not reduce the $ iter(M_L, M_U) $ further. 
On the other side,
more nonzeros can cause the term  $ comp(M_L, M_U) $ increasing, 
and sequentially, the total $ runtime $ increasing.

In order to reduce the runtime,
it should make a balance between the accuracy of SAIT matrices 
and the iteration count.
As the accuracy is mainly decided by 
the threshold $ \tau $ and the number $ p $
in the Algorithm \ref{SAIT_Thr} and \ref{SAIT_Pat}, respectively,
there should be some $ \tau $ or $ p $ 
that make the runtime of an iterative solver with SAIT shortest.
We call such parameters \emph{Optimal parameters}.
The optimal parameters may vary with 
different problems, different methods 
and different computer environments.
We will verify the discussion above using numerical experiments
in section \ref{runtimesait}.

\section{Numerical experiments}

In the numerical experiments,
the main test models are the three-dimensional
Laplace equation and its eigenvalue problem 
with homogeneous Dirichlet boundary condition:
\begin{equation}\label{Lap_00}
	\begin{split}
		-\Delta u &= f \qq\,\,\; \text{in}\,\, \Omega, \qqq
		u = 0 \qq \text{on}\,\, \partial\Omega \\
		\text{and} \qq
		-\Delta u &= \lambda u \qq \text{in}\,\, \Omega, \qqq
		u = 0 \qq \text{on}\,\, \partial\Omega,
	\end{split}
\end{equation}
where $ \Omega = [0,1]^3 $.
We use a $ 102\times 102\times 102 $ uniform mesh
and finite difference method 
to discretize them.
We obtain the corresponding matrix problems
\begin{equation}\label{A_u_f}
	\begin{split}
		A u_h &= f_h\\
		\text{and}  \qq
		A u_h &= \lambda_h u_h.
	\end{split}
\end{equation}
Here, $ A $ is a square matrix with $ 10^6 $ rows/columns
and it is symmetric positive definite (SPD).
We use preconditioned conjugate gradient method (PCG)
and LOBPCG method \cite{knyazev2001toward} to compute
the linear equations and eigenvalue problem in \eqref{A_u_f}, respectively.

We take the level-0 and level-1 ILU factors of the matrix $ A $
as the original triangular systems
and use the SAIT algorithms to approach their inverses.
Then, we take the SAIT matrices as the preconditioners 
in PCG and LOBPCG.
The SAIT matrices are generated by different SAIT algorithms
with different parameters.
We test and compare their effects in preconditioning the iterative solvers.

The code is implemented in Matlab 
and is run on a laptop with an Intel i7-6700HQ CPU with 16 GB RAM
and a GTX970M GPU with 3 GB memory.
We take $ 10^{-10} $ as the uniform stopping criteria
for all the numerical experiments.

\subsection{SAIT\_Thr$(\tau, m)$ and SAIT\_Pat$(p,m)$}

\begin{figure}
	\centering
	\begin{minipage}{6cm}
		\includegraphics[height=6cm,width=6.5cm]{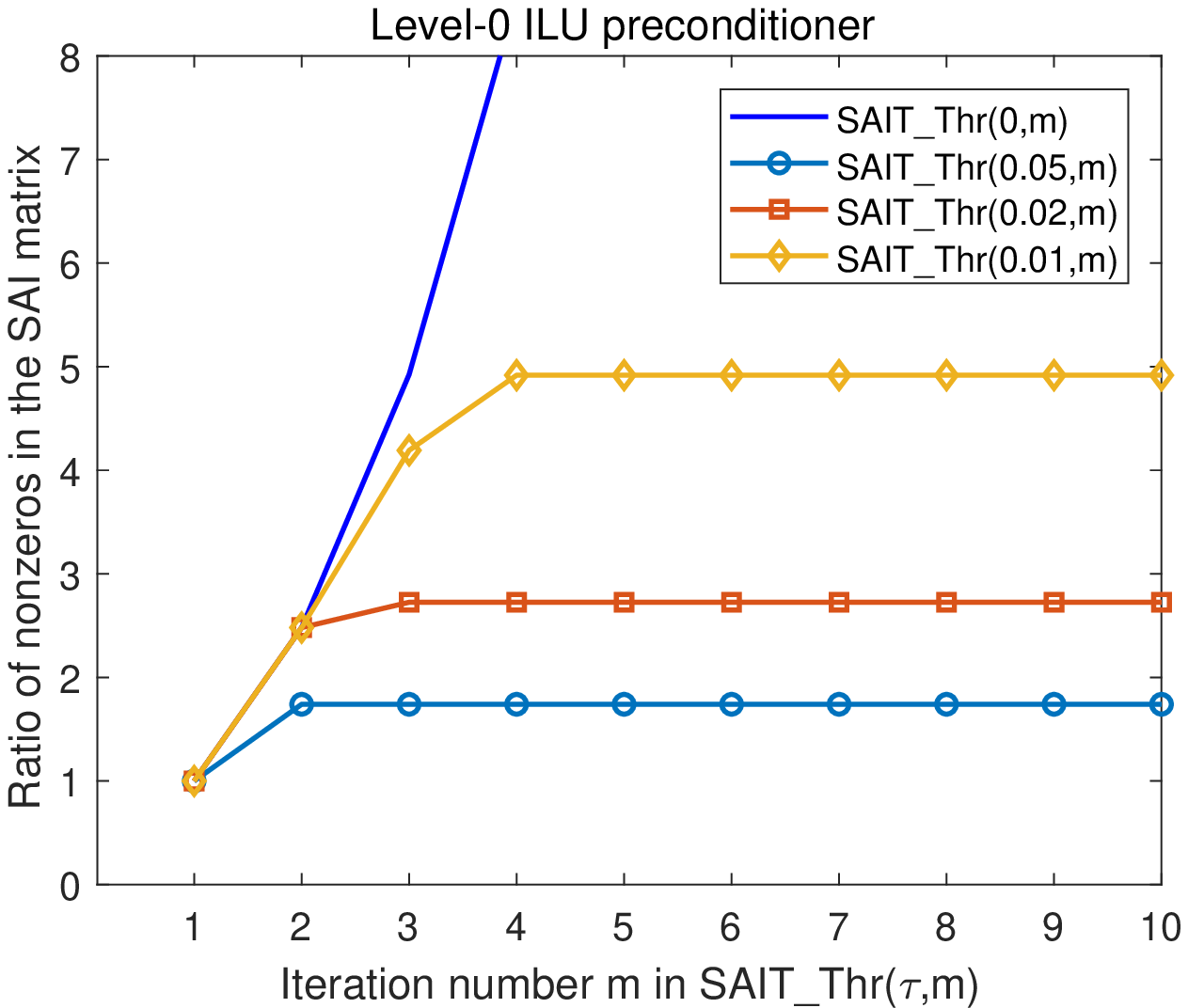}
	\end{minipage}
	\begin{minipage}{6cm}
		\includegraphics[height=6cm,width=6.5cm]{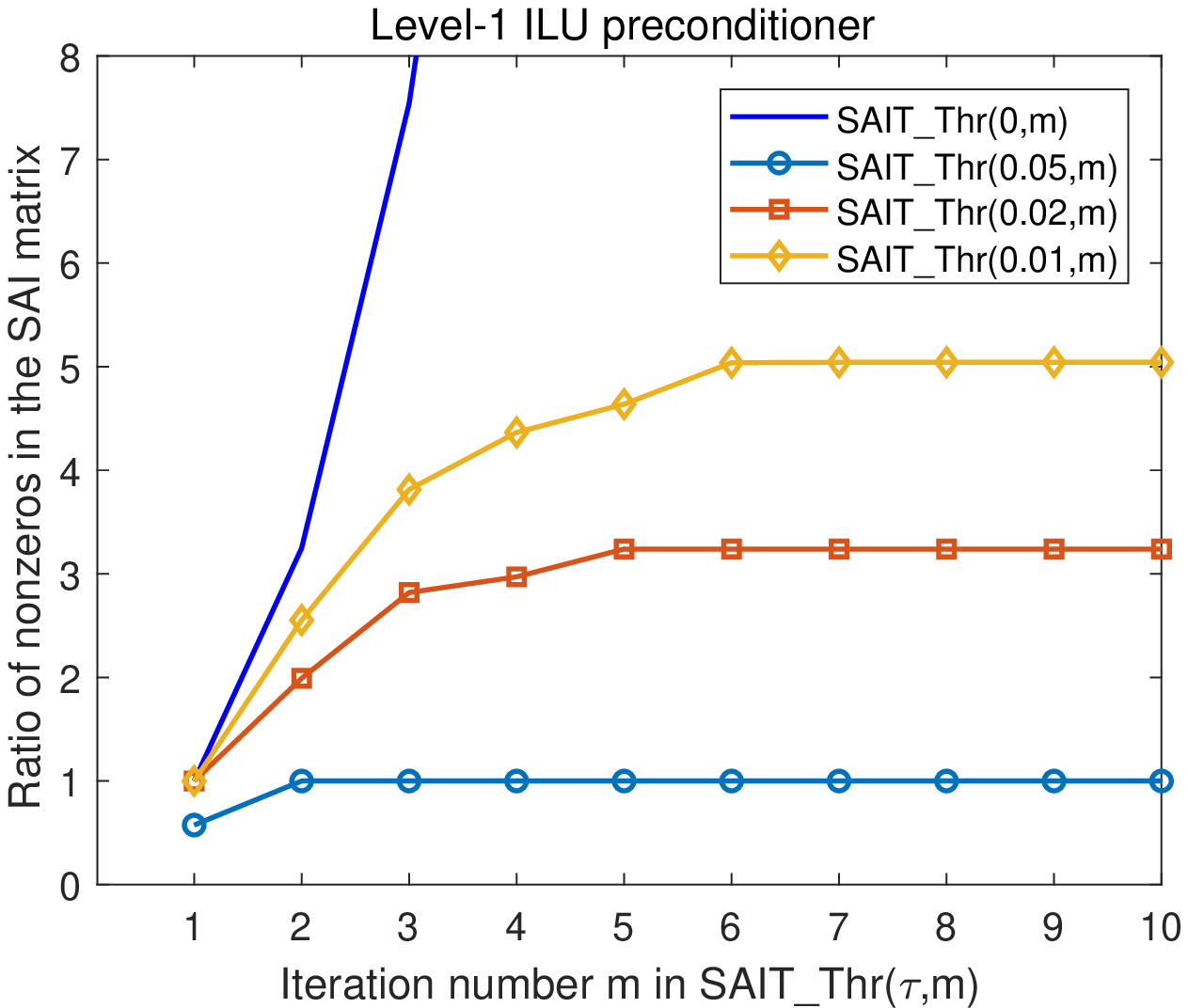}
	\end{minipage}
	\caption{The ratios of the nonzeros in the SAIT matrices generated by
		SAIT\_Thr$(\tau, m)$ with different $ \tau $ and $ m $ for 
		level-0 (left) and level-1 (right) ILU factors.}
	\label{nonzeros_SAI}
\end{figure}

\begin{figure}
	\centering
	\begin{minipage}{6cm}
		\includegraphics[height=6cm,width=6.5cm]{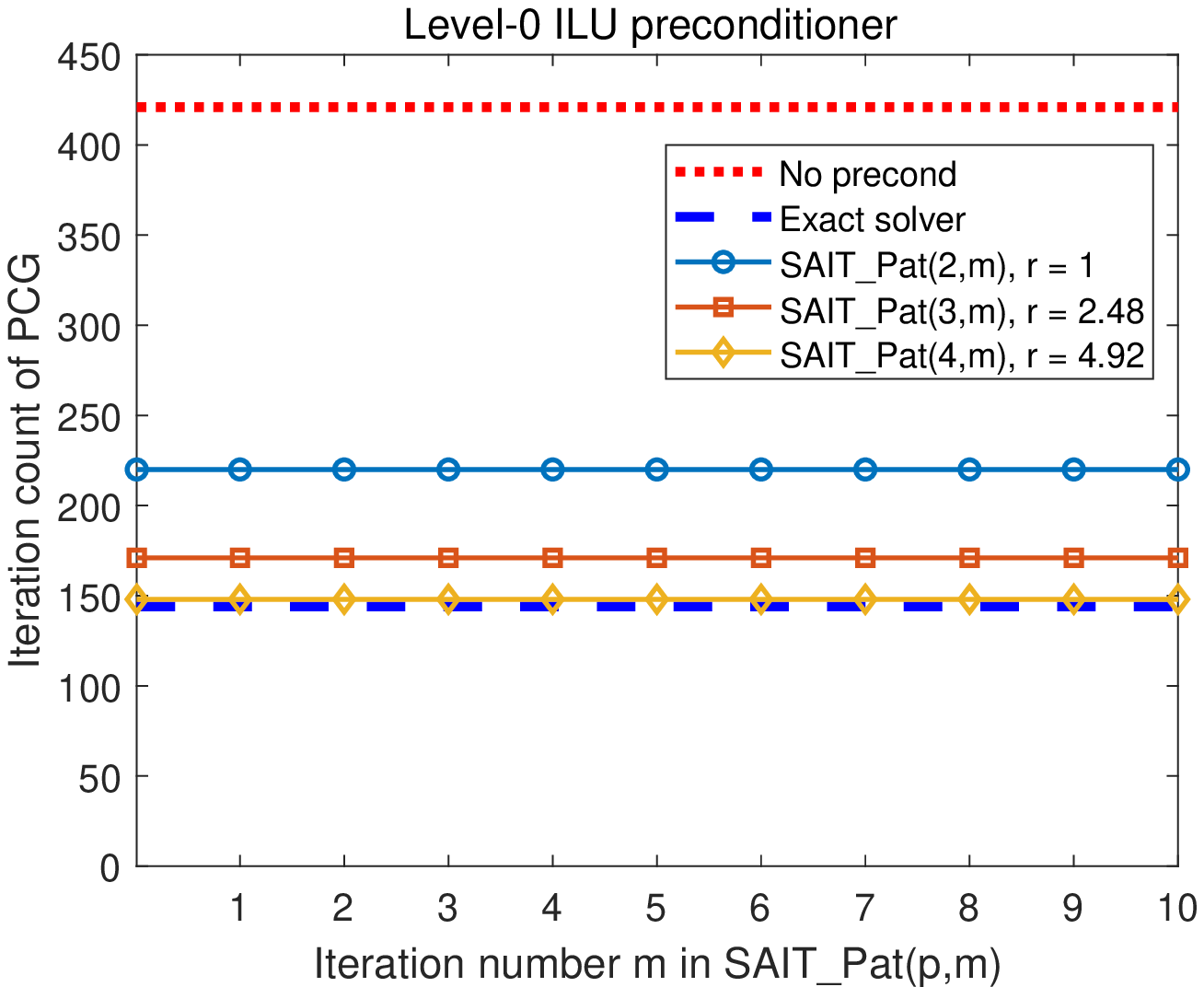}
	\end{minipage}
	\begin{minipage}{6cm}
		\includegraphics[height=6cm,width=6.5cm]{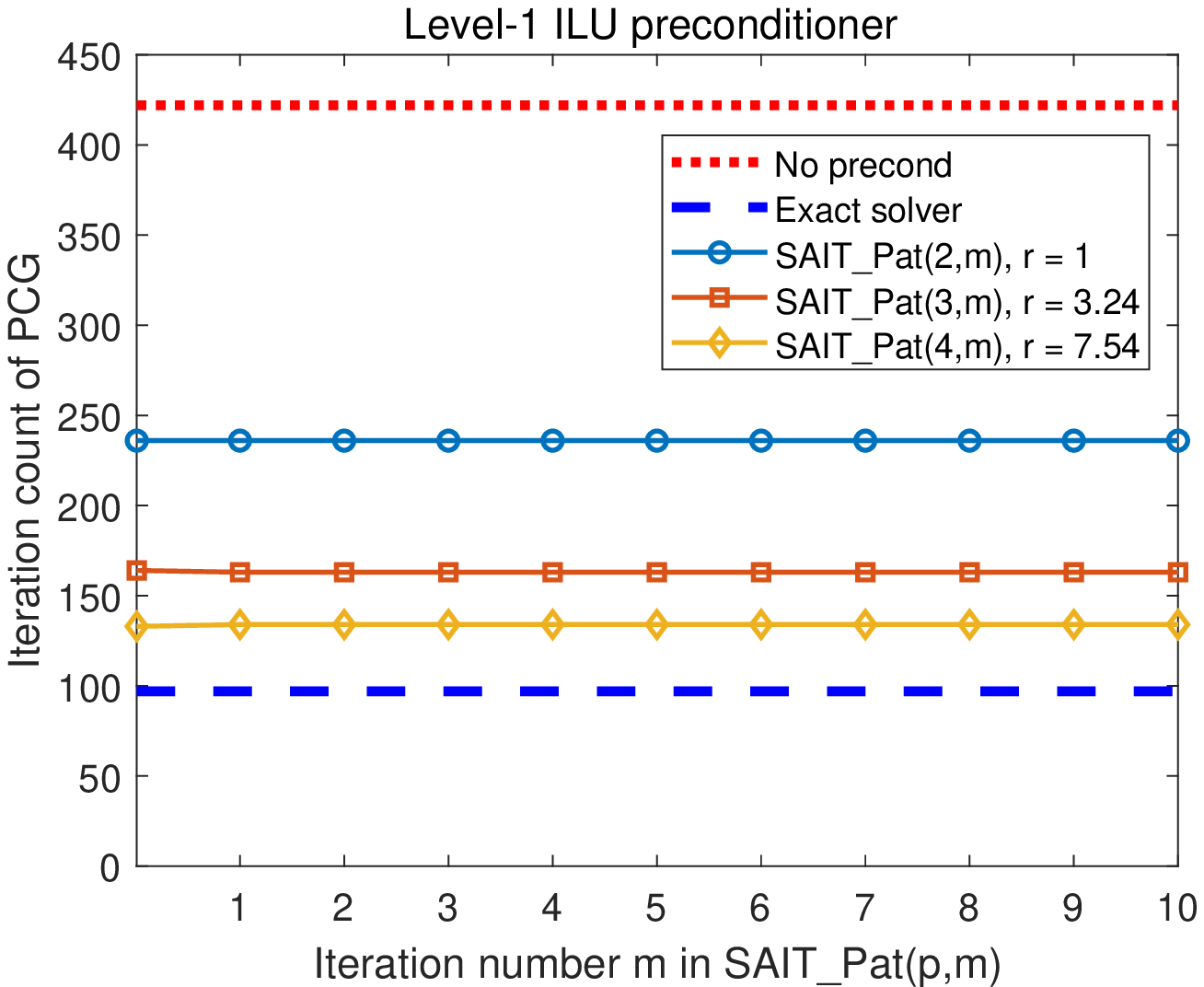}
	\end{minipage}
	\caption{The iteration counts of PCG with SAIT preconditioners generated by
		SAIT\_Pat$(p,m)$ with different $ p $ and $ m $ for 
		level-0 (left) and level-1 (right) ILU factors.}
	\label{PCG_iter_Pat}
\end{figure}

We test the numbers of nonzeros in the SAIT matrices
generated by 
the two dropping strategies,
the SAIT\_Thr$(\tau,m)$ in Algorithm \ref{SAIT_Thr}
and SAIT\_Pat$(p,m)$ in Algorithm \ref{SAIT_Pat}.
Figure \ref{nonzeros_SAI} shows the ratios of nonzeros
in SAIT matrices generated by  SAIT\_Thr$(\tau, m)$
with different threshold $ \tau = 0,\;0.01,\;0.02,\;0.05 $
and increasing iteration number $ m $.
The setting $ \tau =0 $ means that there is no entry dropped in each iteration.
In this case, the numbers of nonzeros increase rapidly.
For other nonzero $ \tau $,
the numbers of nonzeros reach different stable states
after several iterations.
The smaller $ \tau $ results in the more nonzeros in its corresponding SAIT matrix.
For SAIT\_Pat$(p,m)$,
a fixed $ p $ leads to a fixed pattern,
which means that the number of nonzeros of its SAIT matrix is fixed.
We put the ratios of the SAIT matrices generated by this algorithm
in the legends of Figure \ref{PCG_iter_Pat}.

\begin{figure}
	\centering
	\begin{minipage}{6cm}
		\includegraphics[height=6cm,width=6.5cm]{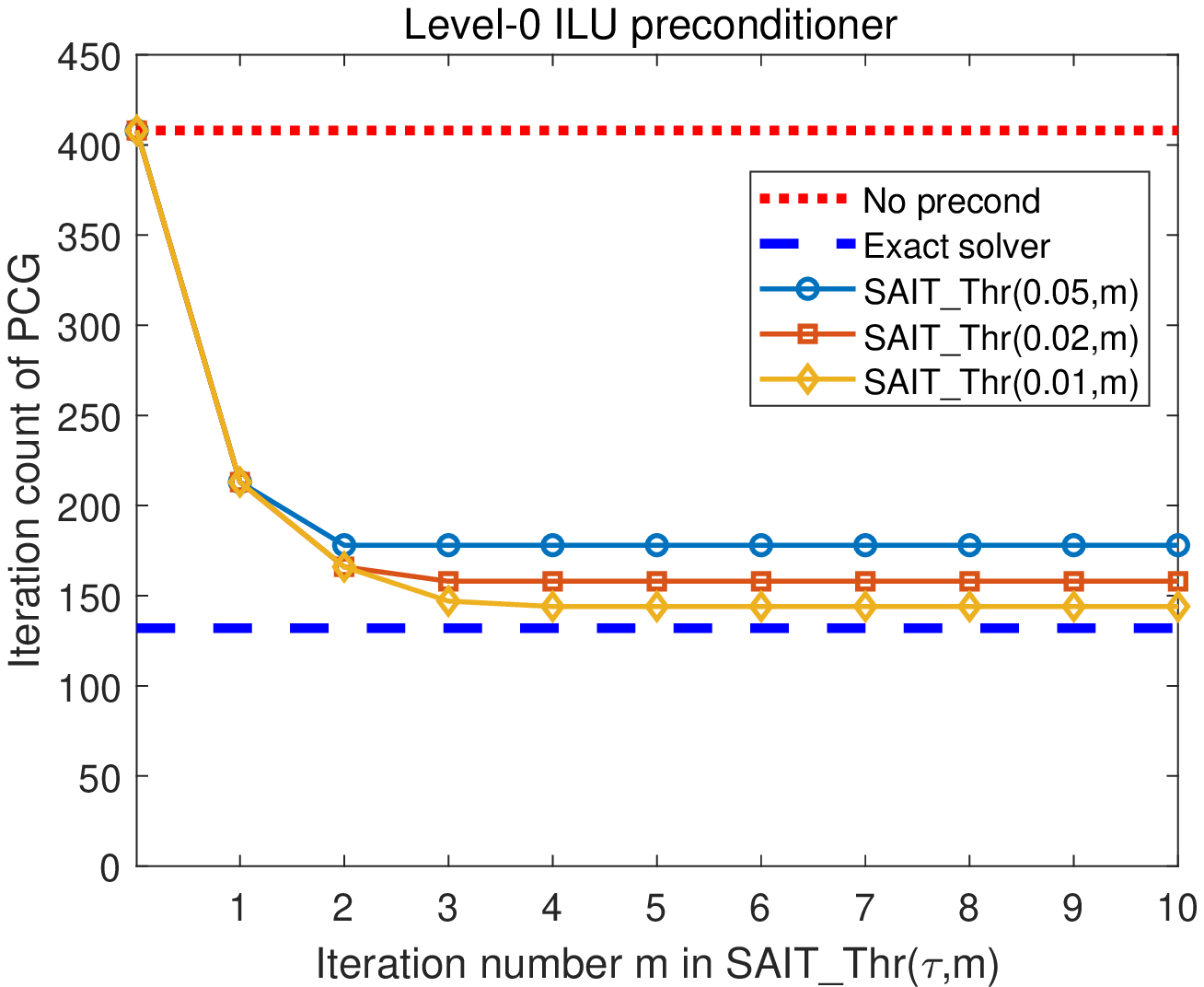}
	\end{minipage}
	\begin{minipage}{6cm}
		\includegraphics[height=6cm,width=6.5cm]{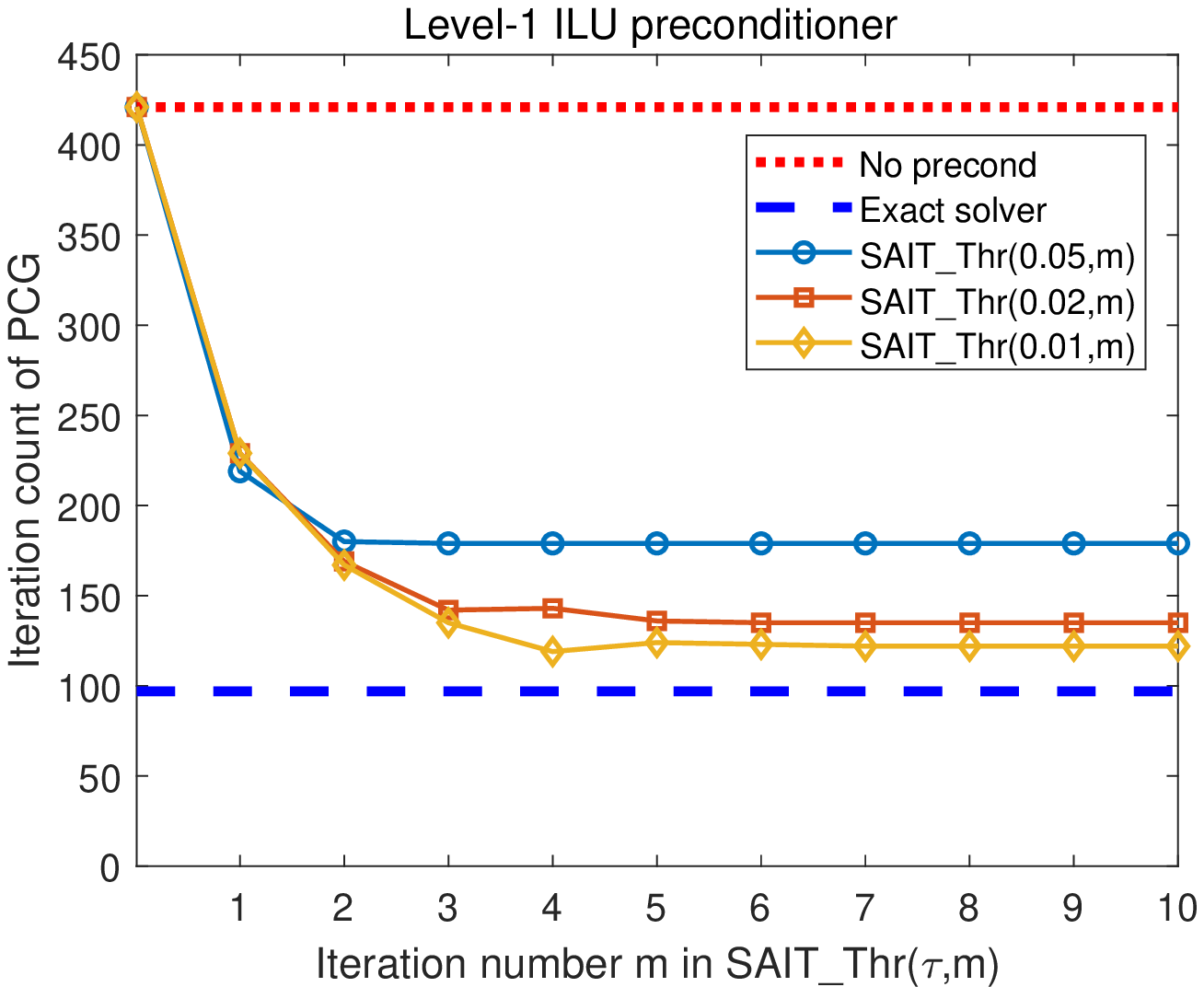}
	\end{minipage}
	\caption{The iteration counts of PCG with SAIT preconditioners generated by
		SAIT\_Thr$(\tau, m)$ with different $ \tau $ and $ m $ for 
		level-0 (left) and level-1 (right) ILU factors. }
	\label{PCG_iter_Thr}
\end{figure}

Next,
we test the effects of SAIT matrices as preconditioners in PCG method.
For the matrices generated by SAIT\_Thr$(\tau, m)$,
the iteration counts of PCG decrease with increasing $ m $
until reaching different stable states, shown in Figure \ref{PCG_iter_Thr}.
When using SAIT\_Pat$(p,m)$,
the iteration counts keep stable with 
$ m $ varying, shown in Figure \ref{PCG_iter_Pat}.
For the both algorithms,
the small threshold $ \tau $ or larger $ p $ (more nonzeros)
means less iteration count.
However, there is no case that 
the SAIT preconditioners can reduced the iteration counts 
less than the exact triangular solvers.

\begin{figure}
	\centering
	\begin{minipage}{6cm}
		\includegraphics[height=6cm,width=6.5cm]{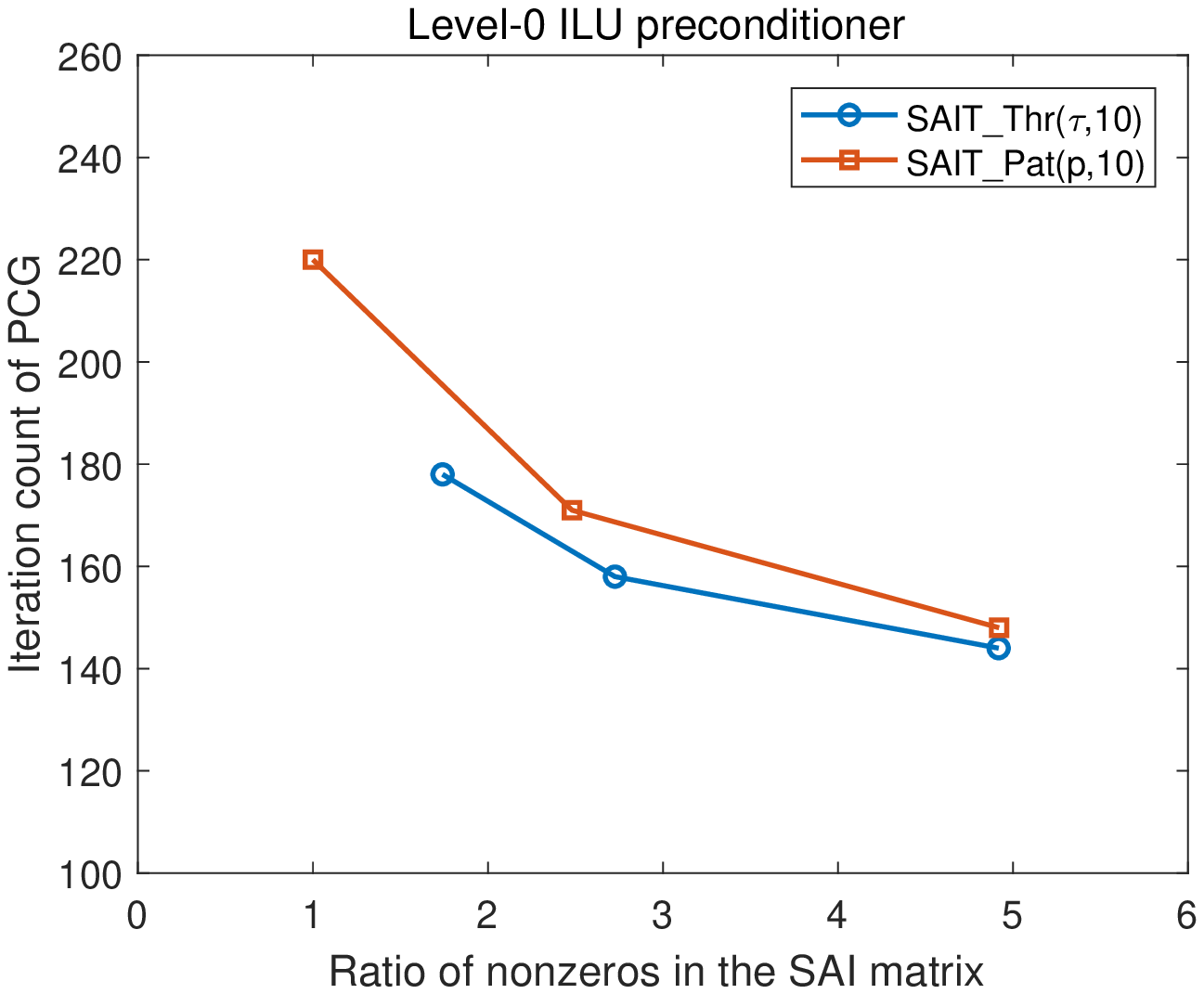}
	\end{minipage}
	\begin{minipage}{6cm}
		\includegraphics[height=6cm,width=6.5cm]{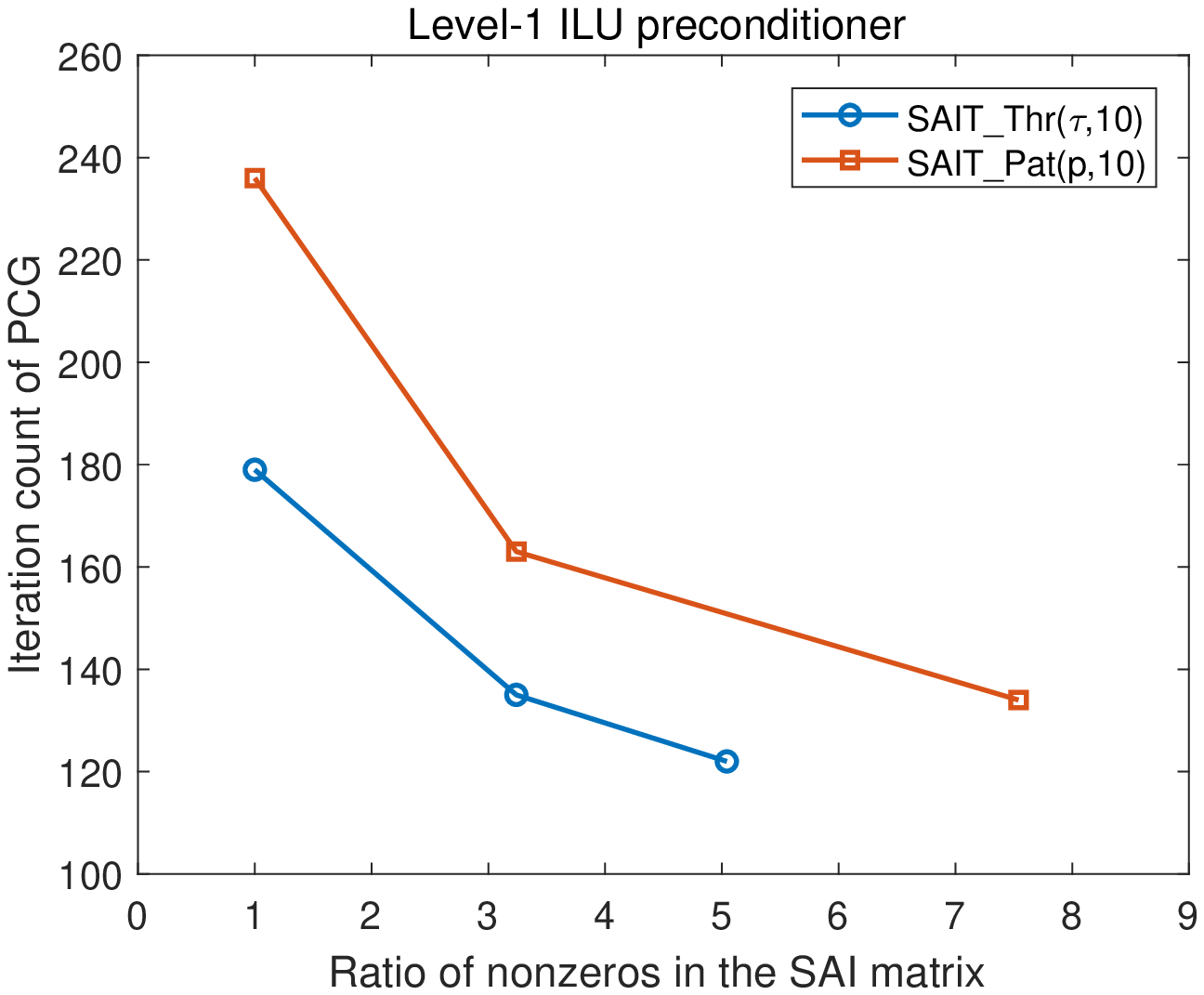}
	\end{minipage}
	\caption{The comparisons of the iteration counts of PCG
		with SAIT preconditioner generated by SAIT\_Thr$(\tau, m)$ and SAIT\_Pat$(p,m)$ 
		for level-0 (left) and level-1 (right) ILU factors. }
	\label{Thr_Pat_compare}
\end{figure}

We also compare the two dropping strategies with a fixed parameter $ m=10 $.
As it is shown in Figure \ref{Thr_Pat_compare},
with the same numbers of nonzeros,
the PCG iteration counts of SAIT\_Thr$(\tau, 10)$
are less than that of SAIT\_Pat$(p,10)$.

\subsection{The balance}

\begin{figure}
	\centering
	\begin{minipage}{6cm}
		\includegraphics[height=6cm,width=6.5cm]{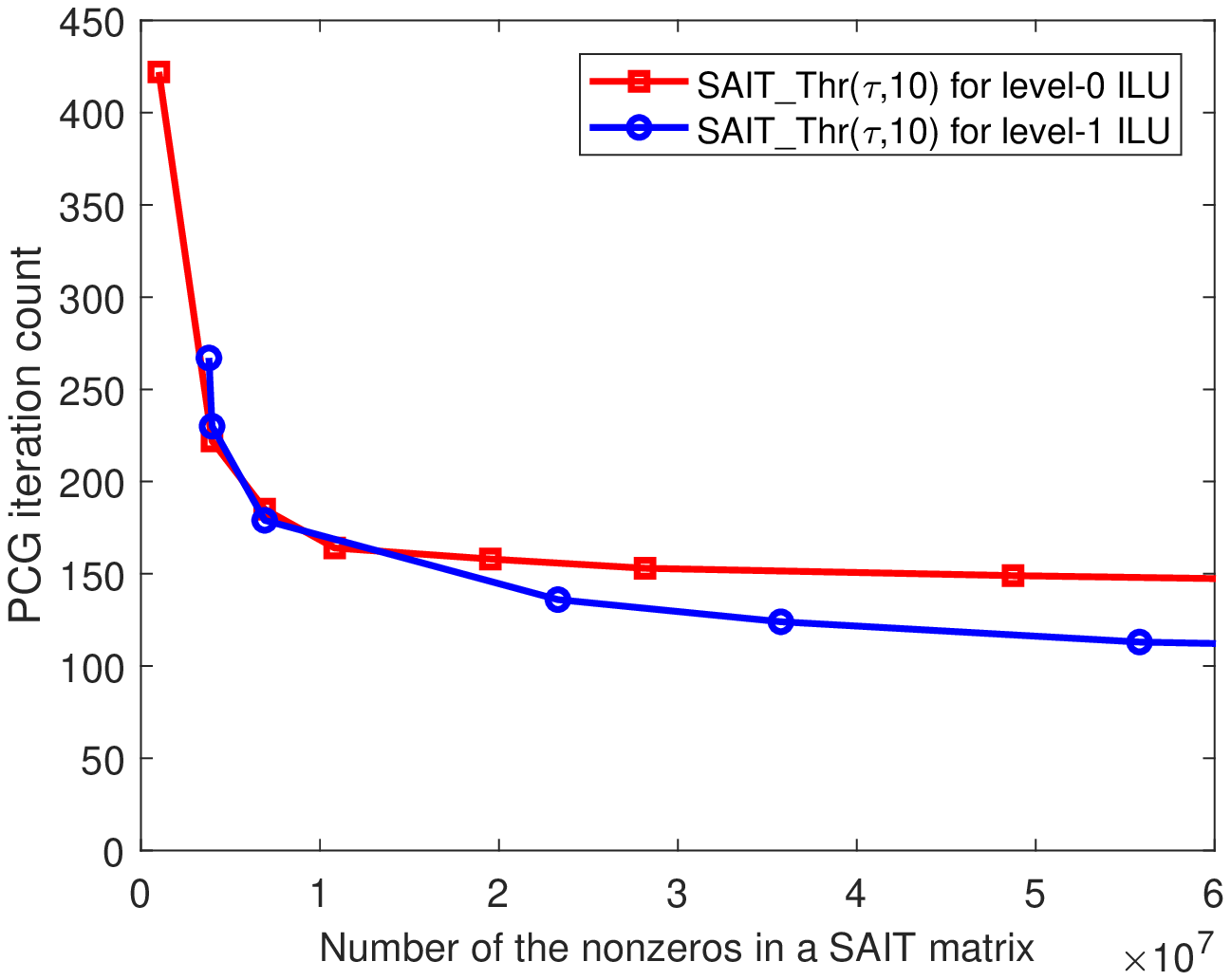}
	\end{minipage}
	\begin{minipage}{6cm}
		\includegraphics[height=6cm,width=6.5cm]{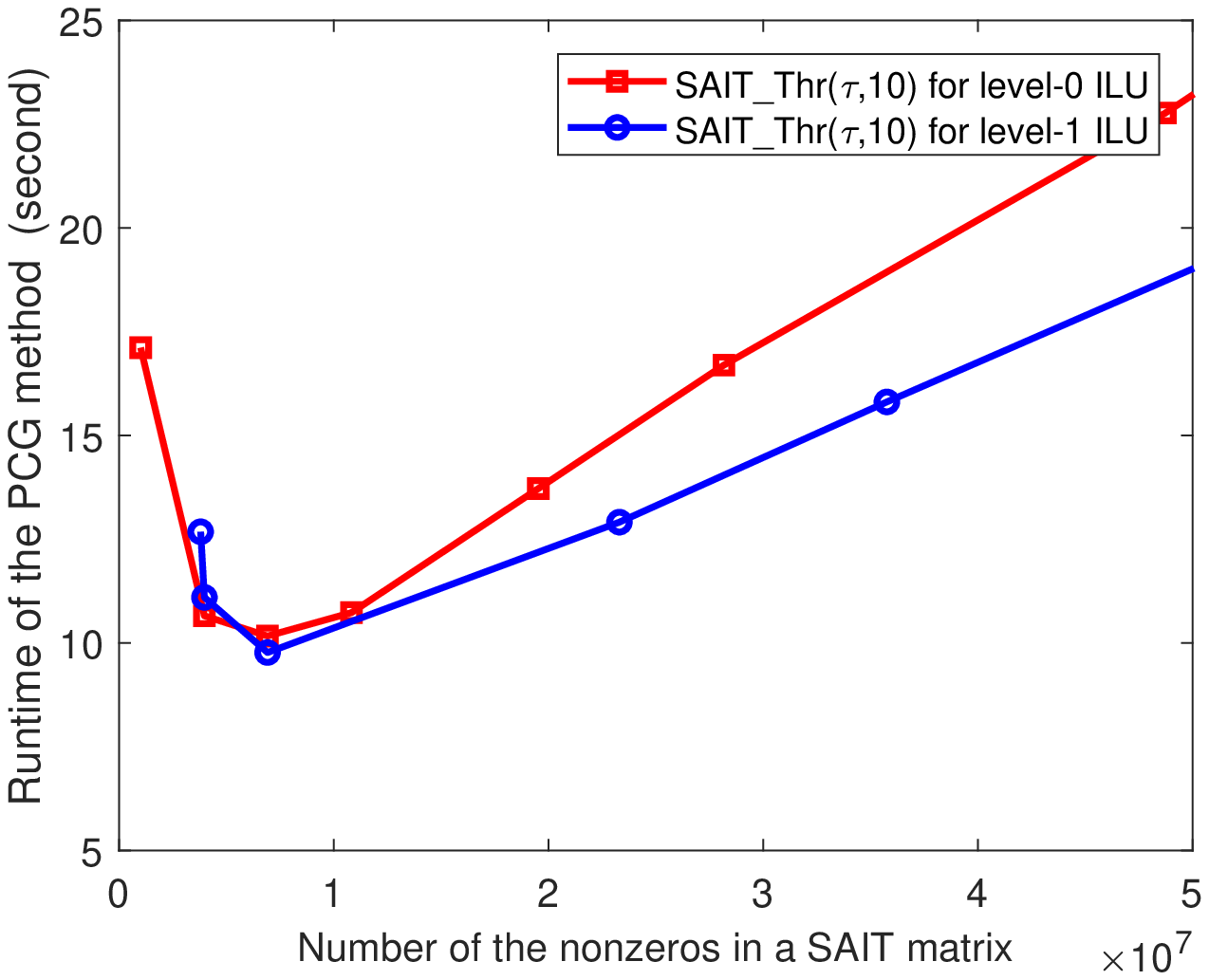}
	\end{minipage}
	\caption{Left: the relation between the PCG iteration counts
		and the numbers of nonzeros in the SAIT preconditioners.
		Right:the relation between the runtime of PCG method
		and the numbers of nonzeros in the SAIT preconditioners.}
	\label{More_SAIT_Thr}
\end{figure}

Table \ref{runtime_SAIT} shows the iteration counts and runtime 
in solving the linear equation in \eqref{A_u_f} 
using PCG with the SAIT preconditioners.
We observe that less iteration count does not means 
shorter runtime.
We test more threshold parameters in SAIT\_Thr$(\tau, m)$,
shown in Figure \ref{More_SAIT_Thr}.
The iteration count does not decrease continuously
with smaller threshold $ \tau $ (more nonzeros) after some point.
Even though we use the exact inverses which are usually full triangular matrices,
the iteration count is equal to the exact triangular solver,
which is not an arbitrarily small number.
However, a SAIT matrix with too many nonzeros can 
increase the cost in computing matrix-vector multiplications.
From the relation between the runtime and numbers of nonzeros,
shown in the right picture of Figure \ref{More_SAIT_Thr},
there should exists an optimal threshold $ \tau $ in respect of the runtime.
Such optimal parameter can be found 
by numerical experiments for each specific problem.

\subsection{Comparison with Jacobi iteration}\label{runtimesait}

Since the SAIT algorithms are designed based on Jacobi iteration method,
if we use the Jacobi method in the preconditioning procedure directly,
it is equivalent to the SAIT matrices generated by SAIT\_Thr$(0, m)$,
which are more accurate than the cases with thresholds $ \tau>0 $.
In addition, 
it requires less memory, 
since a SAIT matrix  usually has more nonzeros than its original matrix.
We list the iteration counts and runtime of the Jacobi method in Table \ref{runtime_Jacobi}.
The PCG iteration counts with Jacobi method can
be almost the same as the exact solver after several iteration.
However, its runtime is much longer than that of SAIT preconditioners
in Table \ref{runtime_SAIT}.
The reason is that
Jacobi method involves several matrix-vector multiplications
and vector-vector additions in each PCG iteration,
while there are only two matrix-vector multiplications 
with the SAIT preconditioners.

\begin{table}
	\centering
	\begin{tabular}{l |c c r |c c r}
		\multirow{ 2}{*}{} &\multicolumn{3}{|c|}{ level-0 ILU } & \multicolumn{3}{|c }{level-1 ILU }\\
		& ratio & iter & time (s) & ratio & iter & time (s)\\
		\hline
		%	Exact solver   &   	-  &  144  &  9.85    &      -  &   98  &  8.31 \\
		SAIT\_Thr(0.05,10)   &     1.74   &  189  &  {\color{red}\emph{\textbf{ 9.95}}} 
		&  1.00  &  184  & {\color{red}\emph{\textbf{ 9.97}}} \\
		SAIT\_Thr(0.02,10)   &     2.73  &  168  &   10.52   &  3.37   &  133  &   12.73  \\
		SAIT\_Thr(0.01,10)   &     4.92  &  {\color{red}\emph{\textbf{  154}}}  &   13.24  &  5.17   & {\color{red}\emph{\textbf{ 121}}}  &   15.48  \\
		\hline
		SAIT\_Pat(1,10)   &       1.00  &  228  &   10.83  &  1.00  &  229  &   12.22  \\
		SAIT\_Pat(2,10)   &    2.48   &  177  &   10.95  &  3.25  &  158  &   15.37  \\
		SAIT\_Pat(3,10)   &    4.92   &{\color{red}\emph{\textbf{  154}}}  &   13.25  &  7.54   &  129  &   22.09 \\
		\hline
	\end{tabular}
	\caption{The ratios the nonzeros, iteration count and runtime in solving the 3D FDM system using PCG method with different SAIT preconditoners. 
		The red italic numbers are the least iteration account or the shortest runtime 
		of the SAIT preconditoners generated
		by level-0 and level-1 ILU factors, respectively.}
	\label{runtime_SAIT}
\end{table}

\begin{table}
	\centering
	\begin{tabular}{c |c c |c c }
		\multirow{ 2}{*}{Jacobi iter} &\multicolumn{2}{|c|}{ level-0 ILU } & \multicolumn{2}{|c }{level-1 ILU }\\
		& iter & time (s) & iter & time (s)\\
		\hline
		1   &     423   &    29.88  & 423   &    31.87  \\
		2   &  229   &    21.03   &  240   &    24.78   \\
		3   &   173   &  {\color{red}\emph{\textbf{19.33}}}   &  169   &    22.10  \\
		4   &    152   &    20.43  &  134   &   {\color{red}\emph{\textbf{ 21.19}}}  \\
		5  &    152   &    23.23  &  116   &    21.48  \\
		6  &    149   &    25.94   &  108   &    23.45   \\
		7  &     147   &    28.19   &  104   &    24.74   \\
		8  &     146   &    30.76  &  101   &    27.21  \\
		9  &     145   &    33.67  &   97   &    28.38  \\
		10  &     145   &    36.60  &   98   &    31.31  \\
		11  &     145   &    39.53  &   98   &    33.90   \\
		12  &     145   &    42.28    &  98   &    36.50   \\
		13  &     145   &    45.66    &  98   &    39.25  \\
		14  &     145   &    48.44   &  97   &    41.43  \\
		15  &     145   &    51.54   &  97   &    47.38   \\
		\hline
	\end{tabular}
	\caption{The iteration counts and runtime 
		in solving the 3D FDM system in \eqref{A_u_f}
		using PCG method with Jacobi iteration method 
		as preconditioner.
		The red italic numbers are the shortest runtime of level-0 and level-1 ILU factors, respectively.}
	\label{runtime_Jacobi}
\end{table}

\subsection{Acceleration by GPU}\label{sec_GPU}

The main computations of the basic Conjecture Gradient method
are matrix-vector multiplications and vector-inner products.
With the SAIT preconditioners,
the preconditioning procedure are two matrix-vector multiplications.
These operations can be highly parallelized.
We use the function $ gpuArray $ in Matlab to input the data into GPU.
We run the solver with the same preconditioners 
in Table \ref{runtime_SAIT} on the GPU.
Though the iteration counts are the similar with the experiments on CPU,
the codes are accelerated several times by GPU,
shown in Table \ref{runtime_SAIT_GPU}.
\begin{table}
	\centering
	\begin{tabular}{l |c c r |c c c}
		\multirow{ 2}{*}{} &\multicolumn{3}{|c|}{ level-0 ILU } & \multicolumn{3}{|c }{level-1 ILU }\\
		& ratio & iter & time (s) & ratio & iter & time (s)\\
		\hline
		%	Exact solver   &   	-  &  144  &  9.85    &      -  &   98  &  8.31 \\
		SAIT\_Thr(0.05,10)   &     1.74   &  189  &    1.34  &  1.00  &  184  &   1.34  \\
		SAIT\_Thr(0.02,10)   &     2.73  &  168  &   1.33   &  3.37   &  133  &   1.59  \\
		SAIT\_Thr(0.01,10)   &     4.92  &  154  &   1.81  &  5.17   &  121   &   1.81 \\
		\hline
		SAIT\_Pat(1,10)   &       1.00  &  228  &   1.47  &  1.00  &  229  &   1.67  \\
		SAIT\_Pat(2,10)   &    2.48   &  177  &   1.36  &  3.25  &  158  &   1.92  \\
		SAIT\_Pat(3,10)   &    4.92   &  154  &   1.80  &  7.54   &  129  &   2.38 \\
		\hline
	\end{tabular}
	\caption{The same codes in Table \ref{runtime_SAIT} are run on GPU.}
	\label{runtime_SAIT_GPU}
\end{table}

\subsection{Preconditioned solver for eigenvalue problems}

\begin{figure}
	\centering
	\begin{minipage}{6cm}
		\includegraphics[height=6cm,width=6.5cm]{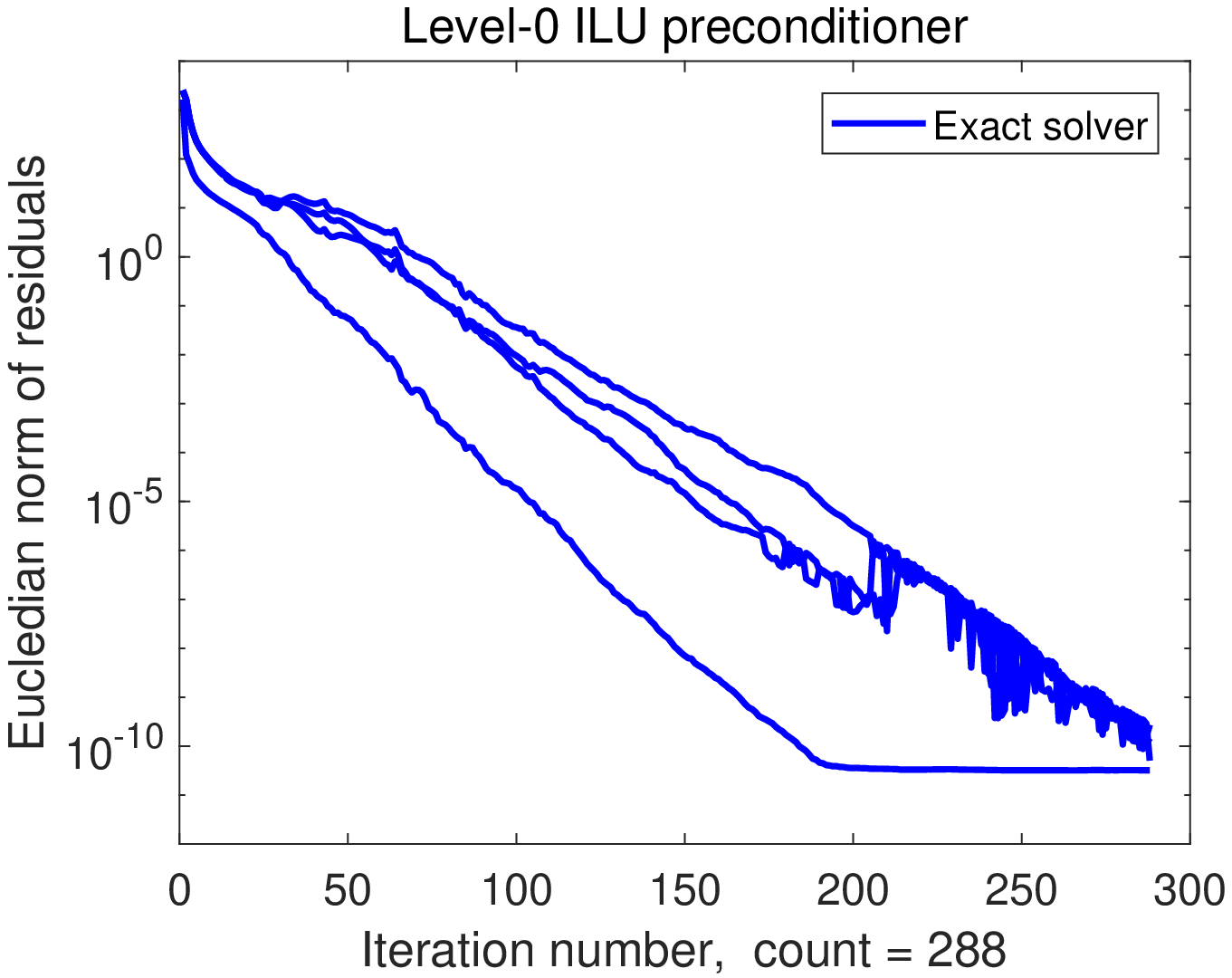}
	\end{minipage}
	\begin{minipage}{6cm}
		\includegraphics[height=6cm,width=6.5cm]{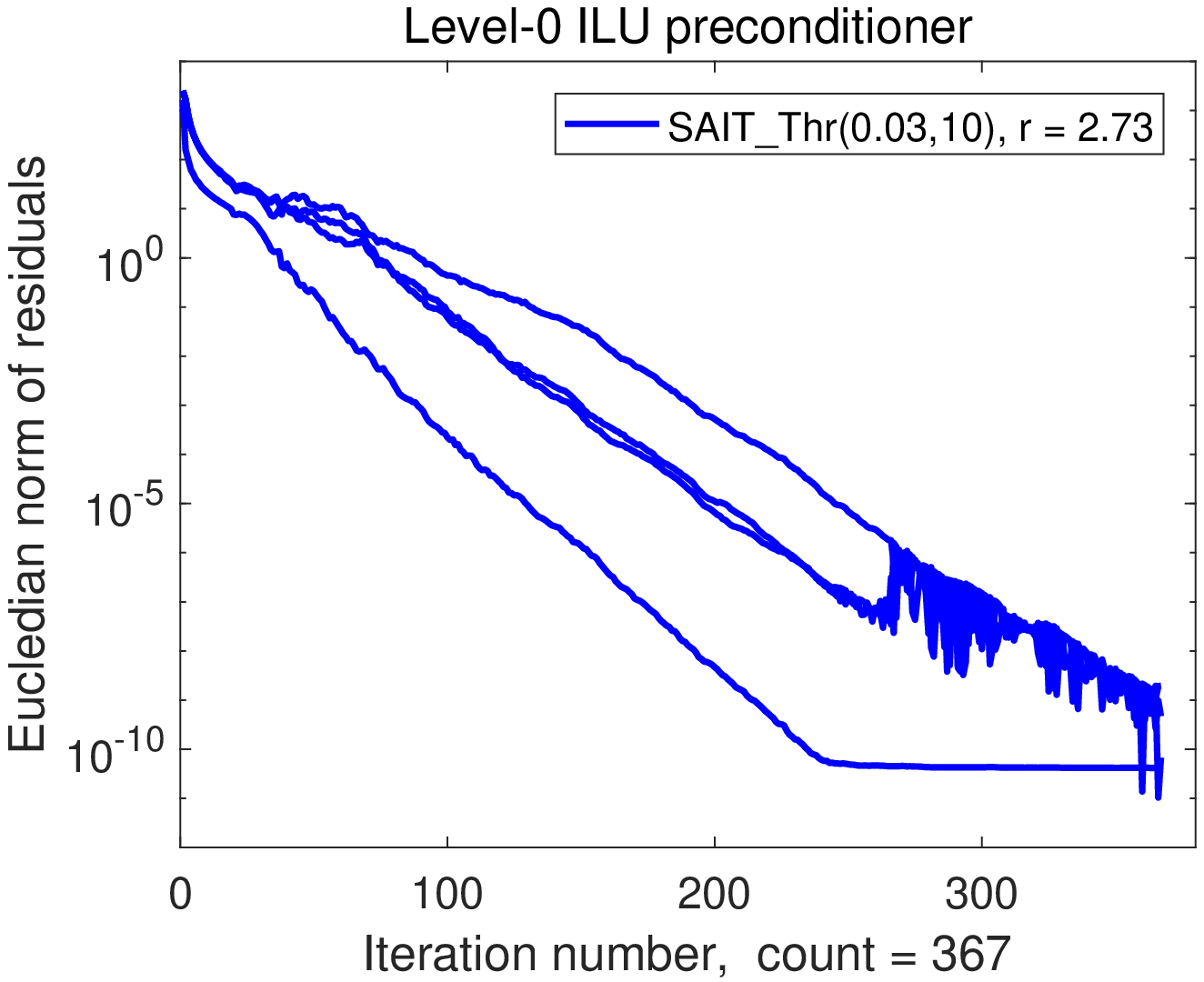}
	\end{minipage}
	\caption{The convergence history in solving Laplace eigenvalue problem
		using LOPBCG method with the exact solver and a SAT preconditioner generate by SAIT\_Thr$(0.03,10)$ for level-0 ILU factors
		in the preconditioning procedure. }
	\label{SAIT_LOPBCG_0}
\end{figure}

\begin{figure}
	\centering
	\begin{minipage}{6cm}
		\includegraphics[height=6cm,width=6.5cm]{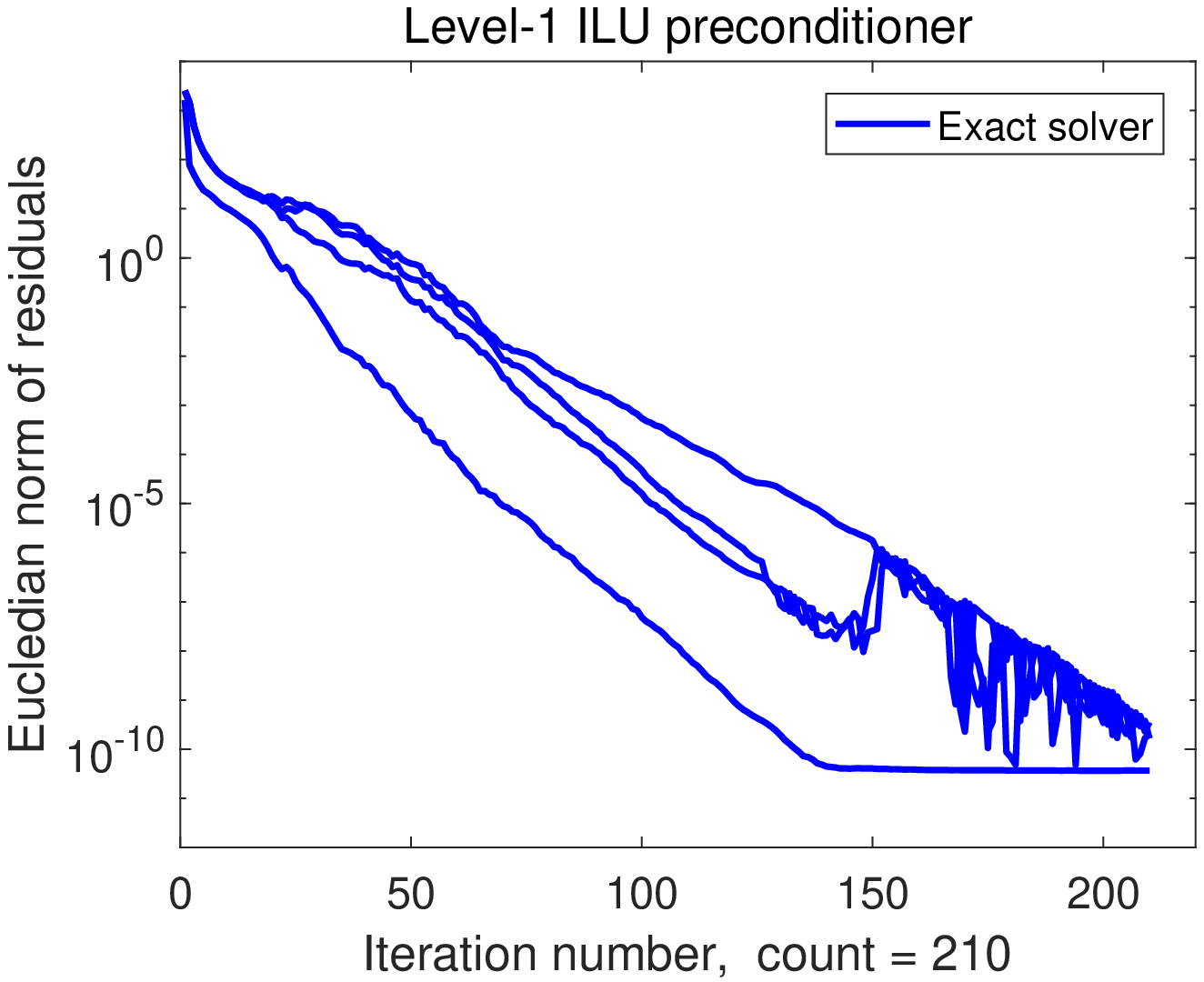}
	\end{minipage}
	\begin{minipage}{6cm}
		\includegraphics[height=6cm,width=6.5cm]{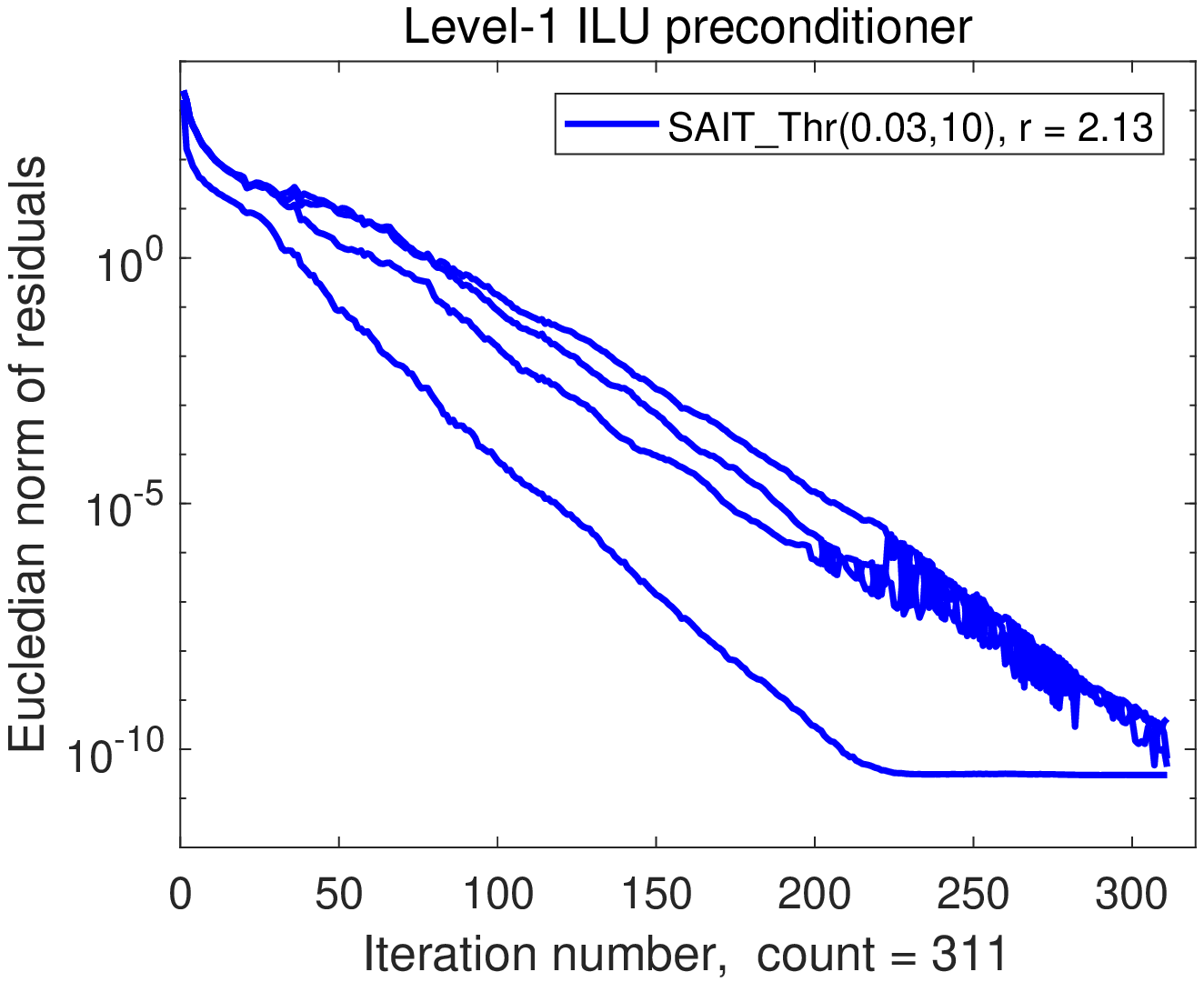}
	\end{minipage}
	\caption{The convergence history in solving Laplace eigenvalue problem
		using LOPBCG method with the exact solver and a SAT preconditioner generated by SAIT\_Thr$(0.03,10)$ 
		for level-1 ILU factors in the preconditioning procedure.}
	\label{SAIT_LOPBCG_1}
\end{figure}

When we use the simultaneous preconditioned method LOBPCG 
to compute the eigenvalue problem in \eqref{A_u_f},
there are many vectors to be dealt with in the preconditioning procedure 
of each iteration.
With SAIT preconditioners,
this procedure can be done through 
two matrix-matrix multiplications.
Here,
we use the LOBPCG method to compute the first 4 eigenvalues
in this problems.
%In order to improve the convergence efficiency,
%we enlarge the subspace in the LOBPCG method 
%by computing one more eigenvector.
We use level-0 and level-1 ILU factorizations,
and use SAIT\_Thr$(0.01, 10)$ to generated the SAIT matrices.
Figure \ref{SAIT_LOPBCG_0} and \ref{SAIT_LOPBCG_1}
show their results, respectively.

\subsection{More examples}
We apply the SAIT algorithms to more examples, 
shown in Table \ref{SPD_matrices_SAIT}.
These examples are from
the University of Florida sparse matrix collection \cite{davis2011university}.
Table \ref{SPD_Thr} are the results using the SAIT\_Thr$(\tau, 10)$ preconditioners.
Here, the ratios of nonzeros in the SAIT matrices generated
by their corresponding thresholds (may be not optimal)
are mainly between 1 and 2.
We also try SAIT\_Pat$(p,10)$ for these examples, shown in Table \ref{SPD_Pat}.
Comparing with the results  in the two tables,
we find that the threshold-based SAIT\_Thr$(\tau, 10)$ preconditioners 
preform  better than the pattern-based SAIT\_Pat$(p,10)$
for these problems.
In Table \ref{SPD_Thr},
we also compare the iteration counts of SAIT preconditioners with the exact solver.
Even though more iterations,
the runtime of SAIT can be reduced by parallelization 
as it is discussed in section \ref{sec_GPU}.
\begin{table}
	\centering
	\begin{tabular}{r| r r }
		& row/column & \ nnz $ \qq $   \\
		\hline
		apache1 & 80800 \ &  542184  \\ 
		apache2 & 715176 \  &  4817870   \\ 
		thermal1 & 82654  \ &  574458   \\ 
		thermal2 & 1228045 \  &  8580313   \\ 
		parabolic\_fem & 525825 \  &  3674625   \\ 
		G3\_circuit & 1585478  \ &  7660826   \\ 
		% McRae/ecology1 & 1000000 &  4996000 &   1 & 1.33 &  1.67 \\ 
		ecology2 & 999999 \  &  4995991  \\ 
		thermomech\_dM & 204316 \  &  1423116    \\ 
		thermomech\_TC & 102158  \ &  711558  \\ 
		%	 offshore & 259789 \  &  4242673 \\ 
		%	 af\_shell3 & 504855 &  17562051 &   1 & 1.42 &  1.97 \\ 
		% Bodendiek/CurlCurl\_3 & 1219574 &  13544618 &   1 & 2.17 &  5.29 \\ 
		\hline
	\end{tabular}
	\caption{The numbers of rows and nonzeros of some SPD matrices in \cite{davis2011university}.}
	\label{SPD_matrices_SAIT}
\end{table}

%\begin{table}
%	\centering
%	\begin{tabular}{r |r| r |r r r | r |r r r }
%		&\multirow{ 2}{*}{no p.c.} &\multicolumn{4}{|c|}{ level-0 ILU } & \multicolumn{4}{|c }{level-1 ILU }\\
%		\cline{3-10}
%		&	& exact & $ \tau $\ \ \ & ratio & SAIT &  exact & $ \tau $ & ratio & SAIT \\
%		\hline
%		apache1 & 3777 & 365 & 0.05 & 1.46 & 439 & 249 & 0.03 & 1.80 & 316 \\
%		apache2 & 5528 & 882 & 0.05 & 1.39 & 1092 & 587 & 0.03 & 1.70  & 797\\
%		thermal1 & 1707 & 651 & 0.05 & 1.43  & 703 & 363 & 0.04 & 1.76 & 435 \\
%		thermal2 & 6626 & 2555 & 0.05 & 1.42 & 2763 & 1401 & 0.04 & 1.72  & 1674 \\
%		parabolic\_fem & 3515 & 1423 & 0.1 \, & 0.96  & 1640 & 845 & 0.04 & 1.32  & 946 \\
%		thermomech\_dM & 89 & 10 & 0.06 & 1.01  & 11 & 6 & 0.02 & 1.02 & 8 \\
%		thermomech\_TC & 89 & 10 & 0.06 & 1.01 & 11 & 6 & 0.02 & 1.02 & 8 \\
%		ecology2 & 7127 & 2123 & 0.08 & 2.00 & 2830 & 1303 & 0.06 & 3.25 & 1942 \\
%		G3\_circuit & 21205 & 1182 & 0.08 & 2.07 & 1582 & 643 & 0.08 & 2.38  & 1174 \\
%		offshore &- & 574 & 0.05 & 1.07 & 690 &- & & & \\
%		
%		\hline
%	\end{tabular}
%	\caption{The iteration counts of PCG method with the exact solver
%	and SAIT matrices generated by SAIT\_Thr$(\tau, m)$ in the preconditioning procedure.}
%	\label{SPD_Thr}
%\end{table}

\begin{table}\small
	\centering
	\begin{tabular}{r |r| r |r r r r | r |r r r r}
		&\multirow{ 2}{*}{no p.c.} &\multicolumn{5}{|c|}{ level-0 ILU } & \multicolumn{5}{|c }{level-1 ILU }\\
		\cline{3-12}
		&	& exact & $ \tau $\ \ \ & ratio & SAIT & &  exact & $ \tau $ & ratio & SAIT & \\
		\hline
		apache1 & 3777 & 365 & 0.05 & 1.46 & 439 & 20.3\%   & 249 & 0.03 & 1.80 & 316& 26.9\%  \\
		apache2 & 5528 & 882 & 0.05 & 1.39 & 1092 & 23.8\%   & 587 & 0.03 & 1.70  & 797& 35.8\% \\
		thermal1 & 1707 & 651 & 0.05 & 1.43  & 703 & 8.0\%   & 363 & 0.04 & 1.76 & 435& 19.8\%   \\
		thermal2 & 6626 & 2555 & 0.05 & 1.42 & 2763 & 8.1\%   & 1401 & 0.04 & 1.72  & 1674& 19.5\%    \\
		parabolic\_fem & 3515 & 1423 & 0.1 \, & 0.96  & 1640 & 15.2\%   & 845 & 0.04 & 1.32  & 946& 12.0\%  \\
		thermomech\_dM & 89 & 10 & 0.06 & 1.01  & 11 & 10.0\%   & 6 & 0.02 & 1.02 & 8&  33.3\% \\
		thermomech\_TC & 89 & 10 & 0.06 & 1.01 & 11 &10.0\%    & 6 & 0.02 & 1.02 & 8&  33.3\% \\
		ecology2 & 7127 & 2123 & 0.08 & 2.00 & 2830 &33.3\%    & 1303 & 0.06 & 3.25 & 1942& 49.0\%   \\
		G3\_circuit & 21205 & 1182 & 0.08 & 2.07 & 1582 & 33.8\%   & 643 & 0.08 & 2.38  & 1174& 82.6\%  \\
		%		offshore &- & 574 & 0.05 & 1.07 & 690& 20.2\% &- &- &- &- &- \\
		\hline
	\end{tabular}
	\caption{The iteration counts of PCG method with the exact solver
		and SAIT preconditioners generated by SAIT\_Thr$(\tau, 10)$.}
	%	The level-1 ILU does not work for the \emph{offshore} problem. }
	\label{SPD_Thr}
\end{table}

\begin{table}
	\centering
	\begin{tabular}{r |r r |r r| r r |r r }
		&\multicolumn{4}{|c|}{ level-0 ILU } & \multicolumn{4}{|c }{level-1 ILU }\\
		\cline{2-9}
		&\multicolumn{2}{|c|}{m = 2}&\multicolumn{2}{|c|}{m = 3}&\multicolumn{2}{|c|}{m = 2}&\multicolumn{2}{|c}{m = 3}\\
		\cline{2-9}
		& ratio & SAIT & ratio & SAIT & ratio & SAIT& ratio & SAIT\\
		\hline
		thermomech\_dM & 1 & 11 & 1.60 &  10 & 1 & 10 & 2.37  & 6  \\
		thermomech\_TC & 1 & 11 & 1.60 &10 & 1 & 10 & 2.37  & 6  \\
		thermal1 & 1 & 847 & 1.85  & 690 & 1 & 743 & 2.49 & 440    \\
		thermal2 & 1 & 3284 & 1.85 & 2674 & 1 & 2868 & 2.48  & 1635    \\
		parabolic\_fem  & 1 & 1678 & 1.56  & 1451 & 1 & 1444 & 2.38 & 893     \\
		apache1   & 1 & 3252 & 2.41 & 2236 & 1 & 4102 & 3.11 & 2570    \\
		apache2   & 1 & 2753 & 2.42 & 1818 & 1 & 3521 & 3.12  & 2500    \\
		G3\_circuit      & 1 & 3993 & 2.04 & 2751 & 1 & 3713 & 2.37  & 2245    \\
		ecology2      & 1 & 3738 & 2.00 & 2799 & 1 & 3765 & 2.25  & 2549    \\
		%		offshore  & 1 &  1389 &     3.96  & 726 &- & - & -&- \\
		\hline
	\end{tabular}
	\caption{The iteration counts of PCG method with the SAIT preconditioners generated by SAIT\_Pat$(p,10)$.}
	\label{SPD_Pat}
\end{table}

\section{Conclusion}

We derive an exact inverse for triangular matrix
trough Jacobi iteration method.
The inverse is a finite series.
We take the truncation of this series as the 
approximate inverse for a triangular matrix.
To make the approximate inverses more sparse,
we propose two dropping strategies.
We apply the SAIT matrices to 
iterative method with ILU preconditioners.
Then the two triangular solvers in ILU preconditioning
are replaced by two matrix-vector multiplications,
which can be fine-grained parallelized.
We present some numerical examples to show 
the effects of SAIT preconditioners.
Though the iteration counts of SAIT matrices
are larger than the exact triangular solvers,
the runtime can be reduced by parallelization.
In this paper,
we only present the basic study on this method.
There are many things about this method are remained
to be studied further,
for example, 
more efficient dropping strategies,
 choices of the parameters
and  parallel implementations.

\begin{appendix}
	
	\section{The Matlab code of SAIT\_Thr$(\tau, m)$} \label{SAIT_Thr_code}

	\lstset{language=Matlab,%
		%basicstyle=\color{red},
		breaklines=true,%
		morekeywords={matlab2tikz},
		keywordstyle=\color{black},%
		morekeywords=[2]{1}, keywordstyle=[2]{\color{black}},
		identifierstyle=\color{black},%
		stringstyle=\color{mylilas},
		commentstyle=\color{mygreen},%
		showstringspaces=false,%without this there will be a symbol in the places where there is a space
		%	numbers=left,%
		%	numberstyle={\tiny \color{black}},% size of the numbers
		%	numbersep=9pt, % this defines how far the numbers are from the text
		emph=[1]{function,for,end,break},emphstyle=[1]\color{blue}, %some words to emphasise
		%emph=[2]{word1,word2}, emphstyle=[2]{style},    
	}
	
	%\begin{mdframed}[tikzsetting={draw=black},skipabove=3pt,]
	%	\begin{Verse}
	%\lstinputlisting{SAIT_Thr.m}
	%	\end{Verse}
	%\end{mdframed}

	%\fbox{
	%\begin{minipage}{20em}
	%\begin{center}
	%\lstinputlisting{SAIT_Thr.m}
	%\end{center}
	%\end{minipage}
	%}

\begin{lstlisting}
function M = SAIT_Thr(T,tau,m)

 I = speye(size(T,1));
Dn = diag(diag(T).^(-1));
T0 = I - Dn*T;
 M = I;

for i = 1 : m
    M = T0*M + I;
    P = (abs(M)>tau);
    M = M.*P;
end

M = M*Dn;

end
\end{lstlisting}

	%\begin{lstlisting}
	%
	%function M = SAIT_Thr(T,tau,m)
	%
	% I = speye(size(T,1));
	%Dn = diag(diag(T).^(-1));
	%T0 = I-Dn*T;
	% M = I;
	%
	%for i = 1 : m
	%    M = T0*M + I;
	%    P = (abs(M)>tau);
	%    M = M.*P;  
	%end
	%
	%M = M*Dn;
	%
	%end
	%
	%\end{lstlisting}

\end{appendix}

\bibliographystyle{plain} 
\bibliography{reference_all} 

\end{document}